**Review of new flow friction equations: Constructing Colebrook's explicit correlations accurately**

**Pavel Praks and Dejan Brkić**

Joint Research Centre (JRC), European Commission, 21027 Ispra, Italy and IT4Innovations, VŠB—Technical University of Ostrava, 708 00 Ostrava, Czechia
Correspondence: pavel.praks@vsb.cz (P.P.) - https://orcid.org/0000-0002-3913-7800; dejanbrkic0611@gmail.com or dejan.brkic@vsb.cz (D.B.) - https://orcid.org/0000-0002-2502-0601
-both authors contributed equally to this study

Abstract: Using only a limited number of computationally expensive functions, we show a way how to construct accurate and computationally efficient approximations of the Colebrook equation for flow friction. The presented approximations are based on the asymptotic series expansion of the Wright ω-function and symbolic regression. The results are verified with 8 million of Quasi-Monte Carlo points covering the domain of interest for engineers. In comparison with the built-in "wrightOmega" feature of Matlab R2016a, the herein introduced related approximations of the Wright ω-function significantly accelerate the computation. With only two logarithms and several basic arithmetic operations used, the presented approximations are not only computationally efficient but also extremely accurate. The maximal relative error of the most promising approximation which is given in the form suitable for engineers' use is limited to 0.0012%, while for a little bit more complex variant is limited to 0.000024%.

**1. Introduction**

The empirical Colebrook equation [1]; Eq. (1)., widely used in hydraulics, implicitly relates in an entangled logarithmic form the known Reynolds number $Re$ and the relative roughness of inner pipe surface $\varepsilon$ with the unknown Darcy flow friction factor $f$. The equation is based on the experiment by Colebrook and White from 1937 [2].

$$\frac{1}{\sqrt{f}} = -2 \cdot log_{10}\left(\frac{2.51}{Re} \cdot \frac{1}{\sqrt{f}} + \frac{\varepsilon}{3.71}\right) \tag{1}$$

To date, the Colebrook equation can be explicitly expressed only in terms of the Lambert *W*-function [3-5], the function introduced by Corless et al. in 1996 [6]. Based on series expansion of its cognate, the Wright ω-function, introduced by Wright [7], explicit approximations of the Colebrook equation of any degree of accuracy can be constructed. With only two computationally expensive logarithmic functions and without exponential nor non-integer power functions, the herein presented approximations which are based on rational polynomial expressions that contain only simple arithmetic operations are also computationally efficient.

Praks and Brkić [8] already developed approximations of the Colebrook equation based on symbolic regression analyses. Symbolic regression is a classic interpretable machine learning method by bridging input data using mathematical expressions composed of some basic functions [9-16]. The study by Praks and Brkić [8] deals with symbolic regression with raw data based on the Colebrook equation solved iteratively and not through the Wright ω-function. Belkić [17] provided a detailed review of applications of transcendental Lambert and Euler functions, which includes mathematics, physics, chemistry, biomedicine, ecology and sociology (hydraulics in [18-20]). Recently, approximations of the Colebrook equation by Brkić and Praks [21,22] based on symbolic regression which goes through the Wright ω-function received a high ranking in large accuracy/speed benchmarks [23,24]. Further comparisons of the available approximations of the Colebrook equation can be seen in [25-35].

Here we will show explicit approximations of the Colebrook equation for flow friction based on the series expansion of $\omega(x) - x$ and symbolic regression. This paper has the following structure. After this

Introduction, Section 2 brings transformation of the Colebrook equation through the Wright ω-function. As the Wright ω-function is computationally demanding, Section 3.1 introduces series expansion of $\omega(x) - x$ using asymptotic expansion and symbolic regression, and then Section 3.2 presents approximations suitable for engineering use, based on series expansions in Section 3.2.1, while on symbolic regression in Section 3.2.2. The novelty of this paper includes also introducing optimal constants, which significantly reduce the relative error of the asymptotic expansion for the Colebrook equation by carefully selected constants, which were estimated by global optimization. Moreover, the novel symbolic regression approximation of the Colebrook equation based on minimization of the relative error is introduced; Section 3.2.2. In comparison with the asymptotic expansion of the same complexity, the novel symbolic-based rational approximation reduces the maximal relative error significantly. Summary related to the presented approximations concerning the accuracy and speed of execution is given in Section 3.3. Finally, Section 4 includes conclusions and future work.

## 2. Transformation through the Wright ω-function

In 1998, Keady [3] was the first to found that the Colebrook equation can be rearranged in the explicit form in terms of the Lambert $W$-function (however Keady's solution has some undesired computational side effects related to a fast-growing term which can cause overflow error [4,5]). As explained by Hayes [36], the Lambert $W$-function, named after the mathematician Johann Heinrich Lambert, who lived in XVIII century, is defined as the inverse of $F(x) = x \cdot e^x$ where $e$ is an exponential function. The Lambert $W$-function has provided to date the unique available closed-form explicit solution of the Colebrook equation in respect to the unknown Darcy's flow friction factor $f$. According to Brkić and Praks [21,22], and based on work by Biberg [37] and by Rollmann and Spindler [38] among others, the Colebrook equation can be exactly explicitly solved by Eq. (2).

$$\left.\begin{array}{c} \frac{1}{\sqrt{f}} = \frac{z}{2.51} \cdot (B + y) \\ z = \frac{2 \cdot 2.51}{\ln(10)} \\ A = \frac{Re}{z} \cdot \frac{\varepsilon}{3.71} \\ B = \ln(Re) - \ln(z) \\ x = A + B \\ y = W(e^x) - x = \omega(x) - x \end{array}\right\} \quad (2)$$

In Eq. (2), the symbol $W$ denotes the Lambert $W$-function and $\omega$ represents the Wright ω-function, where $\omega(x) = W(e^x)$. Further, for the real engineering implementation $y$ denotes an approximation of $\omega(x) - x$, where $\frac{2}{\ln 10} \approx 0.8686$ ; $\frac{2 \cdot 2.51}{\ln 10} \approx 2.18$ ; $2.18 \cdot 3.71 \approx 8.0878$; $A \approx \frac{Re \cdot \varepsilon}{8.0878}$ ; and $B \approx \ln\left(\frac{Re}{2.18}\right) \approx \ln(Re) - 0.7794$.

In this review, we will also test approximations given by the asymptotic expansion of $y \approx \omega(x) - x$, as we confirmed that the accuracy and especially the speed of the algorithm very depends on the used implementation of the Wright ω-function. We will start with approximations based on the asymptotic expansion of the Wright ω-function, and then we will introduce novel approximations of $y$ based on the symbolic regression tool Eureqa [39,40], further noted as $y_{sr}$.

## 3. Methods and discussion

In addition to the already shown Eq. (3) [21] which represents only the first term of series expansion from Eqs. (3-7), we show additional terms in Section 3.1.1 and based on them, related approximations in the form suitable for engineering use in Section 3.2.1. Also, additional strategy on using symbolic regression to construct approximations of $\omega(x) - x$ is given in Section 3.1.2., where their form suitable for engineering use is given in 3.2.2. Besides at the end of Section 3.1.1., we show an efficient strategy on how to reduce the relative error of series expansion-based approximations using carefully selected constants. Only one example of this promising strategy how to construct accurate explicit approximations has been already given in [21], while here we give additional expressions which are to date not only the most accurate but also very simple. These extreme accurate approximations are offered in the form suitable for human eyes (as the number of digits in constants is minimized), but also in forms suitable for coding. Summary of the accuracy and speed of the presented approximations is given in Section 3.3.

### 3.1. Series expansion of $\omega(x) - x$

#### 3.1.1. Assymptotic expansion of $\omega(x) - x$ about infinity

Eqs. (3-7) with $s_i(x)$ terms ($i$=1, 2, …, 5); used for approximating $\omega(x) - x$ by a truncated series about infinity [41], where generally $y_i(x) = \sum_{i=1}^{+\infty} s_i(x)$. The terms $s_1(x), s_2(x), …, s_5(x)$ are given in Eqs. (3-7), see [45]. The corresponding first five approximations of $\omega(x) - x$ are summarized in Eqs. (8-12).

$$s_1(x) = ln(x) \cdot \left(\frac{1}{x} - 1\right) \tag{3}$$

$$s_2(x) = \frac{ln(x)}{2 \cdot x^2} \cdot (ln(x) - 2) \tag{4}$$

$$s_3(x) = \frac{ln(x)}{6 \cdot x^3} \cdot (2 \cdot ln^2(x) - 9 \cdot ln(x) + 6) \tag{5}$$

$$s_4(x) = \frac{ln(x)}{12 \cdot x^4} \cdot (3 \cdot ln^3(x) - 22 \cdot ln^2(x) + 36 \cdot ln(x) - 12) \tag{6}$$

$$s_5(x) = \frac{ln(x)}{60 \cdot x^5} \cdot (12 \cdot ln^4(x) - 125 \cdot ln^3(x) + 350 \cdot ln^2(x) - 300 \cdot ln(x) + 60) \tag{7}$$

As can be seen in Figure 1, the function $y$ needs to be evaluated only in the domain 7.51<$x$<619, which is the one relevant for the Colebrook equation; Eq. (1). It is because the corresponding Reynolds number $Re$ is between 4000 and $10^8$, while the roughness of inner pipe surface $\varepsilon$ is from 0 to 0.05, as it defined in the Moody diagram [42].

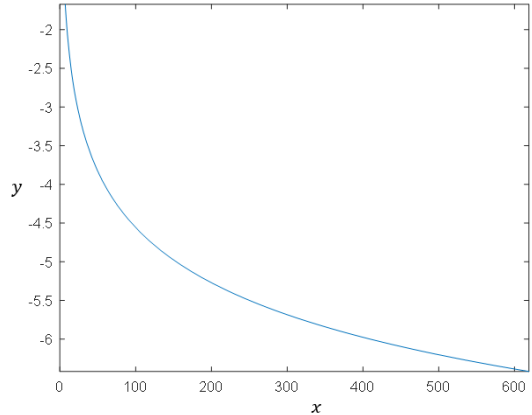

**Figure 1.** The strictly monotonic decreasing function $y = \omega(x) - x$ for the Colebrook equation.

Winning and Coole [43,44] found that logarithmic, exponential, non-integer power terms and the other transcendental functions require extensive additional execution of floating-point operations in the Central Processor Units (CPUs) of modern computers compared with the simple arithmetic operations (+,-,*,/). The logarithmic function is needed only once in Eq. (2); to evaluate $ln(Re)$ and as can be seen from Eqs. (3-7), one additional logarithmic function is needed for evaluation of $ln(x)$ to approximate the $\omega(x)$.

The first five approximations of $\omega(x) - x$; Eqs. (8-12) are based on the series expansion of about infinity by Corless et al. [6,45]. These approximations are expressed in terms of $s_i(x)$ functions from Eqs. (3-7). The coefficients of $\omega(x)$ series expansion about infinity are unsigned Stirling numbers of the first kind [46] as pointed out by Rollmann and Spindler [38].

The symbol $y = wrightOmega(x) - x$ represents the "exact" explicit accurate solution of the Colebrook equation given by Eq. (2), in which $y$ is calculated by the Matlab's built-in command *wrightOmega* (e.g. similar function does not yet exist in MS Excel [47-49] and user-defined Visual Basic for Applications - VBA script should be provided instead [50-55]). However, this *wrightOmega* solution is extremely slow; it requires more than 11 394 seconds, i.e.189.9 minutes(!) to evaluate our benchmark, even the here presented approximations need, depends on the required accuracy, only seconds or even fragments of seconds. In our benchmark study, we measured time in seconds to evaluate the proposed explicit approximations of the Colebrook equation using 8 million of Quasi Monte-Carlo samples from the domain interesting for engineering use. The Monte Carlo method [55-57] for random sampling is based on a broad class of various computational algorithms while for this sampling, the LPTAU51 algorithm was used [58]. We will see from the accuracy analysis and speed tests of Eqs. (8-12) that it makes sense to introduce fast, but accurate approximations of the *wrightOmega* function.

$$y_1 \approx s_1(x) \tag{8}$$
$$y_2 \approx s_1(x) + s_2(x) \tag{9}$$
$$y_3 \approx s_1(x) + s_2(x) + s_3(x) \tag{10}$$
$$y_4 \approx s_1(x) + s_2(x) + s_3(x) + s_4(x) \tag{11}$$
$$y_5 \approx s_1(x) + s_2(x) + s_3(x) + s_4(x) + s_5(x) \tag{12}$$

Eq. (8) is with a maximal relative error of no more than $\{\delta_\%\}_{max}$=0.153% and it needs 0.4 seconds to execute the benchmark sample, Eq. (9) is with $\{\delta_\%\}_{max}$=0.118% and 0.4 seconds, Eq. (10) is with $\{\delta_\%\}_{max}$=0.008% and 1.1 seconds, Eq. (11) with $\{\delta_\%\}_{max}$=0.00249% and 2.4 seconds, while Eq. (12) with $\{\delta_\%\}_{max}$=0.00247% and 4.4 seconds. The maximal error is defined with $\delta_\% = \frac{f - f_i}{f} \cdot 100\%$ where the maximal relative error of the approximated Darcy flow friction factor $f_i$ is calculated with the reference to the accurate $f$ calculated by the Matlab's built-in command *wrightOmega* and Eq. (2). Alternatively, an accurate iterative solution of Eq. (1) is available, as explained in [26,54,59,60].

As can be seen in Eqs. (8-12), the function $y_i$ approximates the function $\omega(x) - x$ by various quality and complexity. The larger $i$ of the function $y_i$ implies a more precise approximation of $\omega(x) - x$, but the computation time increases.

Further, to increase the accuracy of expansions given with Eqs. (8-12), approximations can be written in a form $y + \alpha \approx W(e^x) - x$, where $\alpha$ represents a real constant, which minimizes the relative error of the flow friction factor. The constants were obtained in Matlab by the global optimization solver *Globalsearch*. This optimal constant $\alpha$ is 0.00056 for Eq. (8) gives Eq. (13), -0.0014 for Eq. (9) gives Eq. (14), and -0.000093 for Eq. (10) gives Eq. (15). This strategy also makes the distribution of the error more uniform, and the methodology for the relative error and time measurement is in correlation with Eqs. (8-12).

$$y \approx s_1(x) + 0.00056 \tag{13}$$
$$y \approx s_1(x) + s_2(x) - 0.0014 \tag{14}$$
$$y \approx s_1(x) + s_2(x) + s_3(x) - 0.000093 \tag{15}$$

Eq. (13) is with the maximal relative error of no more than $\{\delta_\%\}_{max}$=0.129% and it needs 0.4 seconds to execute the benchmark sample, Eq. (14) with $\{\delta_\%\}_{max}$=0.0691% and 0.4 seconds, while Eq. (15) with $\{\delta_\%\}_{max}$=0.00527% and 1.1 seconds. From Eq. (13-15) we can observe that the optimal constants tend to zero.

By extending the approach for error minimization by adding a constant as in Eqs. (13-14), a more accurate approximation can be constructed using a function $\xi$ for error minimization (instead of the constant $\alpha$) which is found through symbolic regression [61]. Eq. (16) is based on Eq. (8) where function $\xi$ reduce the relative error significantly (the maximal error is limited to 0.000391% with 1.1 seconds to execute the benchmark).

$$\left. \begin{array}{l} y \approx -ln(x) + \frac{ln(x)}{x} + \xi \\ \xi = \frac{0.3896 \cdot ln(x) \cdot (ln(x) - 1) - 0.9873}{0.8421 \cdot x^2 + 0.01274 \cdot x \cdot (ln(x))^4 + x + 5.882} \end{array} \right\} \tag{16}$$

The error function $\xi$ of Eq. (16) is developed using symbolic regression, while $\omega(x) - x$ will be approximated using the same approach in Section 3.1.2.

### 3.1.2. Symbolic regression-based expansions of $y = \omega(x) - x$ for Colebrook equation

An additional strategy for increasing the accuracy of $y = \omega(x) - x$ for the Colebrook equation is given in Eqs. (17-20). These approximations were found by symbolic regression as already explained in [21] and software Eureqa is used to perform analyses [40] (very recently, also novel artificial intelligence - AI software

AI Feynman: A physics-inspired method for symbolic regression is available [62]). These approximations are optimized for the domain 7.51<$x$<619, which is of the interest for the Colebrook equation.

With the strategy used for Eq. (16), the accuracy of Eq. (19) is significantly improved using a newly constructed error function $\xi_1$ as shown in Eq. (20).

$$y_{SR} \approx -\ln(x) + \frac{1.038 \cdot \ln(x)}{x+0.332} \tag{17}$$

$$y_{SR} \approx -\ln(x) + \frac{1.0119 \cdot \ln(x)}{x} + \frac{\ln(x)-2.3849}{x^2} \tag{18}$$

$$y_{SR} \approx -\ln(x) + \frac{\ln(x)}{x-0.5564 \cdot \ln(x)+1.207} = Y \tag{19}$$

$$y_{SR} \approx Y - \xi_1 = Y - \frac{x \cdot Y^2 + 3.0636 \cdot x \cdot Y + 18.58}{19.5 \cdot (Y^2 \cdot x^2 + x^3) + 169.9 \cdot Y^2 + 1260 \cdot x + 18178} \tag{20}$$

Eq. (17) is with the maximal relative error of no more than $\{\delta_\%\}_{max}$=0.0497% and it needs 0.42 seconds to execute the benchmark sample, Eq. (18) with $\{\delta_\%\}_{max}$ =0.0105% and 0.42 seconds, Eq. (19) with $\{\delta_\%\}_{max}$=0.00229% and 0.42 seconds while Eq. (20) with $\{\delta_\%\}_{max}$=0.000024% and 1.1 seconds.

*3.2. Approximations suitable for engineering use*

3.2.1. Approximations and procedures based on series expansions of $\omega(x) - x$ about infinity

To use the procedure for solving Eq. (2) based on series expansions given with Eqs. (3-7) and Eqs. (8-12), one should follow the algorithm illustrated in Figure 2. It can be modified so $s$ can be given as $s \approx s_1 + s_2$, $s \approx s_1 + s_2 + s_3$, etc. Also, Eqs. (8-12) can be extended to construct $s_6, s_7, s_8$, etc. [45].

Following the algorithm from Figure 2, Eq. (12) can be executed (it introduces a relative error of no more than 0.00247% and that it needs around 4.4 seconds to execute 8 million of Quasi Monte-Carlo tested samples carefully selected from the domain of the Colebrook equation which is of interest for engineering practice).

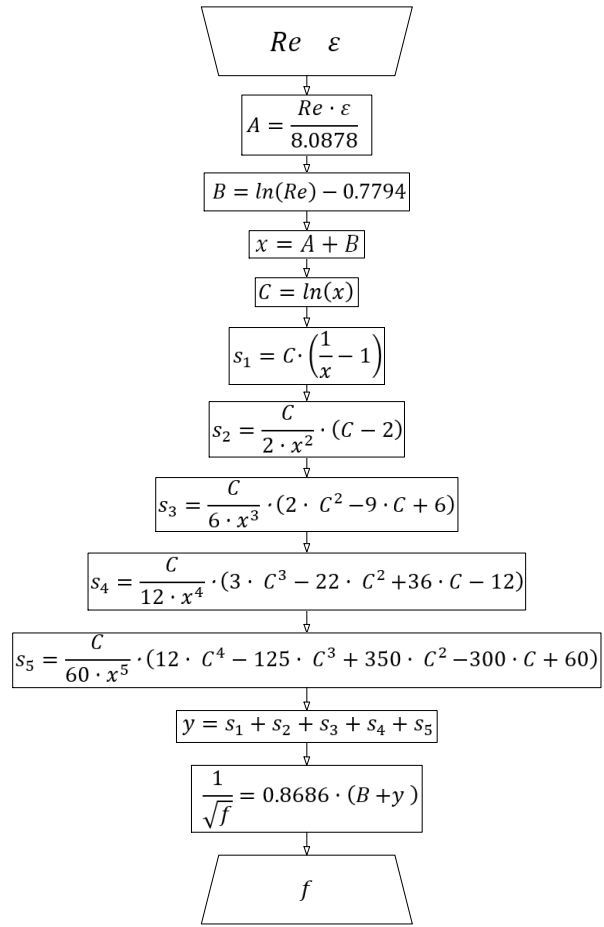

**Figure 2.** Algorithm for the procedure for solving Eq. (2) based on series expansion by Eq. (12).

Besides, if $y \approx s_1$, Eq. (8), the relative error is up to 0.153% with around 0.2 seconds for the execution of 8 million of Quasi Monte-Carlo tested samples. Based on the algorithm from Figure 2 and using only $y \approx s_1$, Eq. (21) can be constructed (this equation is already presented in [21]).

$$\left.\begin{array}{c} \frac{1}{\sqrt{f}} \approx 0.8686 \cdot \left( B - C + \frac{C}{x} \right) \\ A \approx \frac{Re \cdot \varepsilon}{8.0878} \\ B \approx ln(Re) - 0.7794 \\ C = ln(x) \\ x = A + B \end{array}\right\} \qquad (21)$$

Further accurate approximations based on expansions from Eqs. (13-14) are given in Eqs. (22-25). Based on the expansions from Eqs. (13-14), the appropriate explicit approximations to the Colebrook equation can be constructed; Eq. (13) → Eq. (22), Eq. (14) → Eq. (23), and Eq. (15) → Eq. (24).

To construct the approximations from Eqs. (22-24) based on Eqs. (13-14), the notation $C = ln(x)$ should be used to avoid multiple repetitions of the logarithmic function, as $C$ is computed only once.

Comparing Eq. (22) with the maximal relative error of 0.129% and Eq. (21) with 0.153%, it is confirmed that adding the carefully chosen constant $\alpha$, gives more accurate results for the Colebrook equation.

$$\frac{1}{\sqrt{f}} \approx 0.8686 \cdot \left[ B - C + \frac{C}{x} + 0.00056 \right] \qquad (22)$$

$$\frac{1}{\sqrt{f}} \approx 0.8686 \cdot \left[ B - C + \frac{C}{x} + \frac{C}{2 \cdot x^2} \cdot (C - 2) - 0.0014 \right] \qquad (23)$$

$$\frac{1}{\sqrt{f}} \approx 0.8686 \cdot \left[ B - C + \frac{C}{x} + \frac{C}{2 \cdot x^2} \cdot (C - 2) + \frac{C}{6 \cdot x^3} \cdot (2 \cdot C^2 - 9 \cdot C + 6) - 0.000093 \right] \qquad (24)$$

In Eqs. (22-24), $A \approx \frac{Re \cdot \varepsilon}{8.0878}$, $B \approx ln(Re) - 0.7794$, and $C = ln(x)$, $x = A + B$. They are with the relative error of up to $\{\delta_\%\}_{max}$=0.1266%, 0.0691% and 0.00527%, respectively.

Eq. (25), which is based on Eq. (22), brings the relative error of up to 0.000391% thanks to the novel function $\xi$ which reduce the error. Let us compare Eq. (22) and Eq. (25) by terms of accuracy and speed. The relative error of Eq. (25) is reduced significantly by the factor of 0.129%/0.000391% ~ 330. Of course, the computational time of novel Eq. (25) is 1.1 seconds, as it is slower than Eq. (22) by the factor 1.1/0.42~ 2.616, as Eq. (25) is more complex. However, Eq. (25) can be implemented easily in software codes, as $\xi$ is a simple rational function, in which both the numerator and the denominator are polynomials. Here, only the error function $\xi$ is developed using symbolic regression, whereas this technique is more thoroughly evaluated also in Section 3.2.2. The error function $\xi$ is defined in Eq. (16).

$$\left.\begin{array}{c} \frac{1}{\sqrt{f}} \approx 0.86858896 \cdot \left[ B - C + \frac{C}{x} + \varsigma \right] \\ \xi = \frac{0.3896 \cdot C \cdot (C-1) - 0.9873}{0.8421 \cdot x^2 + 0.01274 \cdot x \cdot C^4 + x + 5.882} \\ A \approx \frac{Re \cdot \varepsilon}{8.0884} \\ B \approx ln(Re) - 0.779397 \\ C = ln(x) \\ x = A + B \end{array}\right\} \qquad (25)$$

### 3.2.2. Approximations based on symbolic regression

Four simple but very accurate approximations based on symbolic regression are presented in Eqs. (26-29) where Eqs. (26) and (27) are already given in [21], whereas Eqs. (28) and (29) are novel (they are simple, but very accurate). Although Eureqa allows the only minimization of the absolute error [40], we can use it also for minimizing the relative error as Gholamy and Kreinovich [63] suggested (they showed how to use existing absolute-error-minimizing software to minimize relative errors). Consequently, Eureqa was used for the minimization of the maximal error (worst-case) of $\frac{|y_{sr} - y|}{y}$, where $y_{sr}$ denotes a symbolic regression-based approximation of $y$. Eqs. (26-29) show that this formulation significantly reduces the approximation error.

$$\frac{1}{\sqrt{f}} \approx 0.8686 \cdot \left[ B - C + \frac{1.038 \cdot C}{0.332 + x} \right] \tag{26}$$

$$\frac{1}{\sqrt{f}} \approx 0.8686 \cdot \left[ B - C + \frac{1.0119 \cdot C}{x} + \frac{C - 2.3849}{x^2} \right] \tag{27}$$

$$\frac{1}{\sqrt{f}} \approx 0.8686 \cdot \left[ B - C + \frac{C}{x - 0.5564 \cdot C + 1.207} \right] \tag{28}$$

In Eqs. (26-28), $A \approx \frac{Re \cdot \varepsilon}{8.0884}$, $B \approx ln(Re) - 0.7794$, $C = ln(x)$, and $x = A + B$. They are with the relative error of up to $\{\delta_{\%}\}_{max}$ =0.0497%, 0.0105%, and 0.00337%, respectively (they require only 0.42 seconds to execute the benchmark study).

Comparing Eq. (26) which is already available in [21] with the here introduced Eq. (28), it can be seen that both approximations have similar complexity, but thanks to the minimizing of the relative error [63] through symbolic regression, the maximal relative error of the friction factor was reduced from 0.0497% of Eq. (26) to 0.00337% of Eq. (28), i.e. more than 14 times. Moreover, Eqs. (26-28) are based on Eqs. (17-20), in order to minimize the relative error of the friction factor of the Colebrook equation.

To date, the most accurate approximation of the Colebrook equation is by Vatankhah [64-66] with the relative error of no more than 0.0028% with three logarithmic and one non-integer power functions used (Clamond [67] reported that in a modern computer, evaluation of non-integer power goes also through logarithmic function).

Besides, the accuracy of Eq. (28), can be improved by optimization of coefficients as given in Eq. (29); details of optimization is given in [68,69]. For Eq. (29), the following vector is used [0.8685972, 0.779626, 8.0897, 0.5588, 1.2079] instead of [0.8686, 0.7794, 8.0884, 0.5564, 1.207] for Eq. (28), to decrease the maximal relative error from 0.00337% for Eq. (28) to 0.0012% for Eq. (29). Computation of the benchmark sample of 8 million of Quasi-Monte Carlo points for Eq. (29) takes only 0.388 seconds. However, our novel Eq. (29) with the relative error of no more than 0.0012% is not only two times more accurate, but with two logarithmic functions, it is also simpler (approaches with the only one logarithmic function [70], and even without [71], exist, but they are less accurate; the approach with one logarithm [70] is based on symbolic regression where the second logarithm is replaced using Padé approximants [72], while the algorithm without logarithms [70] is based entirely on rational functions, where both the numerator and the denominator are polynomials). Although the Eq. (28) and Eq. (29) have the same complexity, the optimization of coefficients of Eq. (29) reduced the maximal relative significantly by the factor 0.00337/0.0012~2.8.

$$\frac{1}{\sqrt{f}} \approx 0.8685972 \cdot \left[ B - C + \frac{C}{x - 0.5588 \cdot C + 1.2079} \right] \tag{29}$$

In Eq. (29), $A \approx \frac{Re \cdot \varepsilon}{8.0897}$, $B \approx ln(Re) - 0.779626$, $C = ln(x)$, and $x = A + B$.
Algorithm for solving Eq. (29) is given in Figure 3.

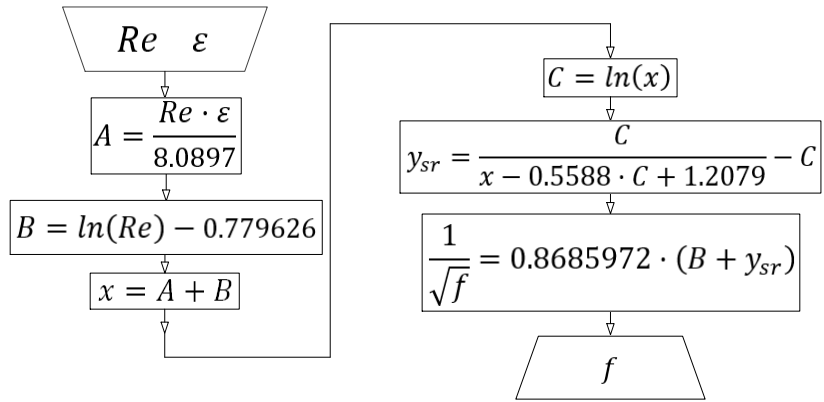

**Figure 3.** Algorithm for solving Eq. (29).

Further, based on Eq. (20), thanks to the function for error minimization $\xi_1$, a novel very accurate approximation given with Eq. (30) can be constructed.

$$\left.\begin{array}{c} \frac{1}{\sqrt{f}} \approx 0.868589 \cdot [B + Y - \xi_1] \\[4pt] \xi_1 = \frac{x \cdot Y^2 + 3.0636 \cdot x \cdot Y + 18.58}{19.5 \cdot (Y^2 \cdot x^2 + x^3) + 169.9 \cdot Y^2 \cdot x + 1260 \cdot x + 18178} \\[4pt] Y = -C + \frac{C}{x - 0.5564 \cdot C + 1.207} \\[4pt] A \approx \frac{Re \cdot \varepsilon}{8.088387} \\[4pt] B \approx ln(Re) - 0.7793975 \\[4pt] C = ln(x) \\[4pt] x = A + B \end{array}\right\} \qquad (30)$$

The computation speed of Eq. (30) is the same as Eq. (15), i.e. 1.1 seconds. The maximal relative error of Eq. (15) is 0.00527%, while the maximal relative error of Eq. 30 is only 0.000024%. Thus, the novel symbolic regression approximation of Eq. (30) reduces the relative error by factor 0.00527%/0.000024%≈219, which confirms the advantages of using symbolic regression for flow friction estimations.

### 3.3. Summary of accuracy and speed tests

The presented approximations of the Colebrook equation which are presented in this article have been compared in Table 1 with respect of accuracy and speed of execution. The column 'time in sec' is given for execution of the Colebrook equation for 8 million samples generated using the Sobol's Quasi Monte Carlo algorithm, as already given through our paper.

Table 1. Comparison of the presented approximations with related error and speed of execution

| $y \approx \omega(x) - x$ | Approximation | Maximal relative error | Time |
|---|---|---|---|
| Exact solution by built-in *wrightOmega* function of Matlab | | | |
| Eq. (2) | | 0% i.e. exact | 11394 seconds ~ 189.9 minutes |
| Asymptotic expansions | | | |
| Eq. (8) | Eq. (21) – Figure 2 | 0.153% | 0.4 seconds |
| Eq. (9) | Figure 2 | 0.118% | 0.4 seconds |
| Eq. (10) | Figure 2 | 0.008% | 1.1 seconds |
| Eq. (11) | Figure 2 | 0.00249% | 2.4 seconds |
| Eq. (12) | Figure 2 | 0.00247% | 4.4 seconds |
| Asymptotic expansions with corrective constant | | | |
| Eq. (13) | Eq. (22) | 0.129% | 0.4 seconds |
| Eq. (14) | Eq. (23) | 0.0691% | 0.4 seconds |
| Eq. (15) | Eq. (24) | 0.00527% | 1.1 seconds |
| Asymptotic expansions with corrective function obtained using symbolic regression | | | |
| Eq. (16) | Eq. (25) | 0.000391% | 1.1 seconds |
| Symbolic regression | | | |
| Eq. (17) | Eq. (26) | 0.0497% | 0.42 seconds |
| Eq. (18) | Eq. (27) | 0.0105% | 0.42 seconds |
| Eq. (19) | Eq. (28) | 0.00337% | 0.42 seconds |
| Symbolic regression with the optimized coefficients | | | |
| Eq. (20); optimized Eq. (19) | Eq. (29) – Figure 3 | 0.0012% | 0.388 seconds |
| Symbolic regression with the optimized coefficients and corrective function | | | |
| Eq. (20) | Eq. (30) | 0.000024% | 1.1 seconds |

Eqs. (8-12) and Eq. (21) are based on asymptotic expansions following algorithm from Figure 2. Eqs. (13-15) and Eqs. (22-24) are based on asymptotic expansions corrected with a carefully selected constant which minimizes the relative error, while Eq. (16) and Eq. (25) are with a carefully selected corrective function. Eqs. (17-19) and Eqs. (26-28) are based on symbolic regression which approximate accurately $\omega(x) - x$; while Eq. (29) is Eq. (28), but with optimised coefficients. Eq. (20) are also based on symbolic regression, but with corrective function. Finally, Eq. (30) is based on accurate approximation of $\omega(x) - x$ found by symbolic regression and also corrected with the carefully selected function for minimization of the maximal relative error (this corrective function is also found using symbolic regression). We can see from Table 1 that it is really worth to construct approximations, as the exact solution of the Colebrook equation; Eq. 2, given by the built-in *wrightOmega* function of Matlab is extremely slow (11 394 seconds, i.e.189.9 minutes(!),)

## 4. Conclusions and future work

In this review, we showed very accurate and fast explicit approximations for the implicitly given Colebrook equation for flow friction. They are based on the Lambert $W$-function and its cognate the Wright $\omega$-function. These explicit approximations are constructed on series expansion and symbolic regression analysis of the Wright $\omega$-function which are much faster compared with the related Matlab built-in function. Moreover, numerical experiments on 8 million of quasi Monte-Carlo pairs show that the approximations based on symbolic regression provide for the fixed complexity more accurate flow friction approximations than the general approximations based on the series expansion. The reason is that the here presented symbolic regression approximations are optimized for the Colebrook equation, whereas the series expansion is more general.

The proposed approximations should be compared with very recent comparative studies of the available approximations of the Colebrook equation [23-26]. Such comparisons will confirm that the approximations proposed in this review are both accurate and computationally efficient.

In terms of the balance among simplicity for human eyes, accuracy and computational simplicity, the novel Eq. (29) with a relative error of up to 0.0012% is highly recommended (a simple algorithm shown in Figure 3 should be followed for engineering use). For a very accurate software implementation, its more complex variant given by Eq. (30) should be followed, as the maximal relative error is only 0.000024%.

**Contributors:** Both authors contributed equally to this study. Pavel Praks worked mostly on the mathematical background of the proposed methods and procedures, while Dejan Brkić put them in a form suitable for engineering use. Pavel Praks was responsible for experiments with symbolic regression tools and global optimization. Dejan Brkić made the draft of this paper.

**Funding:** This work has been partially supported by the Technology Agency of the Czech Republic under the project "National Centre for Energy" TN01000007,"Energy System for Grids" TK02030039 and by the Ministry of Education, Youth and Sports of the Czech Republic through the National Programme of Sustainability (NPS II) project "IT4Innovations excellence in science-LQ1602".

**Acknowledgments:** We especially thank Marcelo Masera for his careful scientific supervision and kind approvals. We acknowledge support from the European Commission, Joint Research Centre (JRC), Directorate C: Energy, Transport, and Climate, Unit C3: Energy Security, Distribution, and Markets.

## Notations

The following symbols are used in this paper:

| | |
|---|---|
| $A$ | variable that depends on $Re$ and $\varepsilon$ (dimensionless); $A = \frac{Re}{z} \cdot \frac{\varepsilon}{3.71}$ |
| $B$ | variable that depends on $Re$ and $z$ (dimensionless); $B = ln(Re) - ln(z) \approx ln(Re) - 0.7794$ |
| $C$ | variable that depends on $x$ (dimensionless); $C = ln(x)$ |
| $f$ | Darcy (Moody) flow friction factor (dimensionless) – main output parameter |
| $s$ | approximation of $y$ by a truncated series about infinity; $y_i = \sum_{i=1}^{+\infty} s_i$ |
| $Re$ | Reynolds number (dimensionless) – main input parameter |
| $x$ | variable that depends on $A$ and $B$ (dimensionless); $x = A + B$ |
| $y$ | approximation of $\omega(x) - x$, i.e. of $W(e^x) - x$ |
| $y_{sr}$ | approximation of $\omega(x) - x$, i.e. of $W(e^x) - x$ found by symbolic regression |
| $Y$ | approximation of $\omega(x) - x$ from Eq. (19) and used in Eq. (20) |
| $z$ | constant (dimensionless); $z = \frac{2 \cdot 2.51}{ln(10)}$ |
| $\alpha$ | constant (dimensionless); half of the smallest term in the asymptotic series |
| $\varepsilon$ | Relative roughness of inner pipe surface (dimensionless) – main input parameter |
| $\delta_\%$ | Relative error (%); $\delta_\% = \frac{f - f_i}{f} \cdot 100\%$ |
| $\xi$ and $\xi_1$ | error-functions developed using symbolic regression |
| $e$ | exponential function |
| $log_{10}$ | logarithm with base 10 |
| $ln$ | natural logarithm |
| $W$ | Lambert $W$-function |
| $\omega$ | Wright $\omega$-function |
| $i$ | Counter |
| *wrightOmega* | Matlab's built-in function |
| *Globalsearch* | Matlab's built-in solver |

## References

[1]. Colebrook C.F. Turbulent flow in pipes with particular reference to the transition region between the smooth and rough pipe laws. J. Inst. Civ. Eng., 11:133–156, 1939. doi:10.1680/ijoti.1939.13150

[2]. Colebrook C., White C. Experiments with fluid friction in roughened pipes. Proc. R. Soc. Lond. Ser. A Math. Phys. Sci., 161:367–381, 1937. doi:10.1098/rspa.1937.0150

[3]. Keady G. Colebrook-White formula for pipe flows. J. Hydraul. Eng., 124:96–97, 1998. doi:10.1061/(ASCE)0733-9429(1998)124:1(96)

[4]. Sonnad J.R., Goudar C.T. Constraints for using Lambert W function-based explicit Colebrook–White equation. J. Hydraul. Eng., 130:929-931, 2004. doi:10.1061/(ASCE)0733-9429(2004)130:9(929)

[5]. Brkić D. Comparison of the Lambert W-function based solutions to the Colebrook equation. Eng. Comput., 29:617-630, 2012. doi:10.1108/02644401211246337

[6]. Corless R.M., Gonnet G.H., Hare D.E., Jeffrey D.J., Knuth D.E. On the Lambert W function. Adv. Comput. Math., 5:329–359, 1996. doi:10.1007/BF02124750

[7]. Wright E.M. Solution of the equation $z \cdot e^z = a$. Bull. Am. Math Soc., 65:89–93, 1959. doi:10.1090/S0002-9904-1959-10290-1

[8]. Praks P., Brkić D. Symbolic regression-based genetic approximations of the Colebrook equation for flow friction. Water, 10:1175, 2018. doi:10.3390/w10091175

[9]. Niazkar M., Talebbeydokhti N., Afzali S.-H. Bridge backwater estimation: A comparison between artificial intelligence models and explicit equations. Sci. Iran., doi:10.24200/SCI.2020.51432.2175

[10]. Filippov A.E., Gorb S.N. Methods of the pattern formation in numerical modeling of biological problems. Facta Universitatis, Series: Mech. Eng., 17:217-242, 2019. doi:10.22190/FUME190227027F

[11]. Gorkemli B., Karaboga D. A quick semantic artificial bee colony programming (qsABCP) for symbolic regression. Inf. Sci. (N.Y.), 502:346–362, 2019. doi:10.1016/j.ins.2019.06.052

[12]. Niazkar M. Revisiting the estimation of Colebrook friction factor: A comparison between artificial intelligence models and CW based explicit equations. KSCE J. Civ. Eng., 23:4311–4326, 2019. doi:10.1007/s12205-019-2217-1

[13]. Diyaley S., Chakraborty S. Optimization of multi-pass face milling parameters using metaheuristic algorithms. Facta Universitatis, Series: Mech. Eng., 17:365-383, 2019. doi:10.22190/FUME190605043D

[14]. Brkić D., Ćojbašić Ž. Evolutionary optimization of Colebrook's turbulent flow friction approximations. Fluids, 2:15, 2017. doi:10.3390/fluids2020015

[15]. Brkić D., Ćojbašić Ž. Intelligent flow friction estimation. Comput. Intel. Neurosc., 2016:5242596, 2016. doi:10.1155/2016/5242596

[16]. Ćojbašić Ž., Brkić D. Very accurate explicit approximations for calculation of the Colebrook friction factor. Int. J. Mech. Sci., 67:10-13, 2013. doi:10.1016/j.ijmecsci.2012.11.017

[17]. Belkić Dž. All the trinomial roots, their powers and logarithms from the Lambert series, Bell polynomials and Fox–Wright function: illustration for genome multiplicity in survival of irradiated cells. J. Math. Chem., 57:59–106, 2019. doi:10.1007/s10910-018-0985-3

[18]. Brkić, D. Lambert W function in hydraulic problems. Math. Balk. New Ser., 26:285–292, 2012, http://www.math.bas.bg/infres/MathBalk/MB-26/MB-26-285-292.pdf

[19]. Brkić D. New explicit correlations for turbulent flow friction factor. Nucl. Eng. Des., 241:4055-4059, 2011, doi:10.1016/j.nucengdes.2011.07.042

[20]. Brkić, D. W solutions of the CW equation for flow friction. Appl. Math. Lett., 24:1379-1383, 2011, doi:10.1016/j.aml.2011.03.014

[21]. Brkić D., Praks P. Accurate and efficient explicit approximations of the Colebrook flow friction equation based on the Wright ω-function. Mathematics, 7:34, 2019. doi:10.3390/math7010034

[22]. Brkić D., Praks P. Accurate and efficient explicit approximations of the Colebrook flow friction equation based on the Wright ω-function: Reply to Discussion. Mathematics, 7:410, 2019. doi:10.3390/math7050410

[23]. Zeyu Z., Junrui C., Zhanbin L., Zengguang X., Peng L. Approximations of the Darcy–Weisbach friction factor in a vertical pipe with full flow regime. Water Sci. Technol.: Water Supply, 20, 2020. doi:10.2166/ws.2020.048

[24]. Muzzo L.E., Pinho D., Lima L.E., Ribeiro L.F. Accuracy/Speed analysis of pipe friction factor correlations. In International Congress on Engineering and Sustainability in the XXI Century—INCREaSE 2019, Section: Water for Ecosystem and Society, Faro, Portugal, 9–11 October 2019; Monteiro, J., Silva, A.J., Mortal, A., Aníbal, J., da Silva, M., Oliveira, M., Eds.; Springer Nature Switzerland AG 2020: Cham, Switzerland, 2019; pp. 664–679, doi:10.1007/978-3-030-30938-1_51

[25]. Brkić D. Review of explicit approximations to the Colebrook relation for flow friction. J. Pet. Sci. Eng., 77:34–48, 2011. doi:10.1016/j.petrol.2011.02.006

[26]. Brkić D. Determining friction factors in turbulent pipe flow. Chem. Eng. N.Y. 119:34-39, 2012.

[27]. Anaya-Durand A.I., Cauich-Segovia G.I., Funabazama-Bárcenas O., Gracia-Medrano-Bravo V.A., Evaluación de ecuaciones de factor de fricción explícito para tuberías (Evaluation of explicit friction factor equations for pipes) /in Spanish/. Educ. Quím., 25:128-134, 2014. doi:10.1016/S0187-893X(14)70535-X

[28]. Silva A.L.M.F., Mesquita A.L.A., Felipe A.M.P.F., Souza, J.A.D.S. Simulação do fator de atrito para o escoamento confinado de caulim com diferentes teores de sólidos pelo modelo de SWAMEE–JAIN (Simulation of the friction factor for the confined flow of kaolin with different solid contents using the SWAMEE – JAIN model) /in Portuguese/. Rev. Mater., 22, 2017. doi:10.1590/s1517-707620170002.0146

[29]. Pimenta B.D., Robaina A.D., Peiter M.X., Mezzomo W., Kirchner J.H., Ben L.H. Performance of explicit approximations of the coefficient of head loss for pressurized conduits. Rev. Bras. Eng. Agríc. Ambient., 22:301-307, 2018. doi:10.1590/1807-1929/agriambi.v22n5p301-307

[30]. Santos-Ruiz I., Bermúdez J.R., López-Estrada F.R., Puig V., Torres, L. Estimación experimental de la rugosidad y del factor de fricción en una tubería /in Spanish/. In Proceedings of the Memorias del Congreso Nacional de Control Automático, San Luis Potosí, Mexico, 10–12 October 2018; Available online: http://amca.mx/RevistaDigital/cnca2018/articulos/VieAT3-03.pdf

[31]. Olivares A., Guerra R., Alfaro M., Notte-Cuello E., Puentes L. Evaluación experimental de correlaciones para el cálculo del factor de fricción para flujo turbulento en tuberías cilíndricas (Experimental evaluation of correlations used to calculate friction factor for turbulent flow in cylindrical pipes) /in Spanish/. Rev. Int. Métodos Numér. Cálc. Diseño Ing., 35:15, 2019. doi:10.23967/j.rimni.2019.01.001

[32]. Alfaro-Guerra M., Guerra-Rojas R., Olivares-Gallardo A. Evaluación experimental de la solución analítica exacta de la ecuación de Colebrook-White (Experimental evaluation of exact analytical solution of the Colebrook-White Equation) /in Spanish/. Ingeniería Investigación y Tecnología, 20:1–11. 2019. doi:10.22201/fi.25940732e.2019.20n2.021

[33]. Minhoni R.T.D.A., Pereira, F.F., da Silva, T.B., Castro, E.R., Saad, J.C. The performance of explicit formulas for determining the Darcy-Weisbach friction factor. Eng. Agric., 40:258-265, 2020. doi:10.1590/1809-4430-eng.agric.v40n2p258-265/2020

[34]. Pérez-Pupo J.R., Navarro-Ojeda M.N., Pérez-Guerrero J.N., Batista-Zaldívar M.A. On the explicit expressions for the determination of the friction factor in turbulent regime. Revista Mexicana de Ing. Quim., 19:313–334, 2020. doi:10.24275/rmiq/Fen497

[35]. Niazkar M., Talebbeydokhti, N. Comparison of explicit relations for calculating Colebrook friction factor in pipe network analysis using h-based methods. Iran. J. Sci. Technol., Trans. Civ. Eng., 44:231–249, 2020. doi:10.1007/s40996-019-00343-2

[36]. Hayes B. Why W? Am. Sci., 93:104–108, 2005.

[37]. Biberg, D. Fast and accurate approximations for the Colebrook equation. J. Fluids Eng., 139:031401, 2017. doi:10.1115/1.4034950

[38]. Rollmann P., Spindler K. Explicit representation of the implicit Colebrook–White equation. Case Stud. Therm. Eng., 5:41–47, 2015. doi:10.1016/j.csite.2014.12.001

[39]. Schmidt M., Lipson H. Distilling free-form natural laws from experimental data. Science, 324:81–85, 2009. doi:10.1126/science.1165893

[40]. Dubčáková R. Eureqa: Software review. Genet. Program. Evolvable Mach., 12:173–178, 2011. doi:10.1007/s10710-010-9124-z

[41]. Horchler A.D. Complex double-precision evaluation of the Wright Omega function. 2017. https://github.com/horchler/wrightOmegaq

[42]. Moody L.F. Friction factors for pipe flow. T. ASME, 66:671–684.

[43]. Winning H.K., Coole T. Explicit friction factor accuracy and computational efficiency for turbulent flow in pipes. Flow Turbul. Combust., 90:1–27, 2013. doi:10.1007/s10494-012-9419-7

[44]. Winning H.K., Coole T. Improved method of determining friction factor in pipes. Int. J. Numer. Methods Heat Fluid Flow, 25:941–949, 2015. doi:10.1108/HFF-06-2014-0173

[45]. Corless R.M., Jeffrey D.J. The Wright ω Function. In: Calmet J., Benhamou B., Caprotti O., Henocque L., Sorge V. (eds) Artificial Intelligence, Automated Reasoning, and Symbolic Computation. AISC 2002, Calculemus 2002. Lecture Notes in Computer Science, vol 2385. Springer, Berlin, Heidelberg, 2002. doi:10.1007/3-540-45470-5_10

[46]. On-Line Encyclopedia of Integer Sequences; Sequence A008275 for Stirling numbers of first kind. http://oeis.org/A008275

[47]. Brkić D. Spreadsheet-based pipe networks analysis for teaching and learning purpose. Spreadsheets Educ. (EJSiE), 9:4646, 2016. https://sie.scholasticahq.com/article/4646

[48]. Niazkar M., Afzali S.H. Analysis of water distribution networks using MATLAB and Excel spreadsheet: h-based methods. Comput. Appl. Eng. Educ., 25:129-141, 2017. doi:10.1002/cae.21786

[49]. Niazkar M., Afzali S.H. Analysis of water distribution networks using MATLAB and Excel spreadsheet: Q-based methods. Comput. Appl. Eng. Educ., 25:277-289, 2017. doi:10.1002/cae.21796

[50]. Demir S., Karadeniz A., Demir N.M., Duman S. Excel VBA-based solution to pipe flow measurement problem. Spreadsheets Educ. (EJSiE), 10:4671, 2018. https://sie.scholasticahq.com/article/4671.pdf

[51]. Niazkar M., Afzali S.H. Streamline performance of Excel in stepwise implementation of numerical solutions. Comput. Appl. Eng. Educ., 24:555–5664, 2016. doi:10.1002/cae.21731

[52]. Niazkar M., Talebbeydokhti N., Afzali, S.H. Relationship between Hazen-William coefficient and Colebrook-White friction factor: Application in water network analysis. Eur. Water, 58:513-520, 2017. https://www.ewra.net/ew/pdf/EW_2017_58_74.pdf

[53]. Niazkar M., Talebbeydokhti N., Afzali S.H. Novel grain and form roughness estimator scheme incorporating artificial intelligence models. Water Resour. Manag., 33:757-773, 2019. doi:10.1007/s11269-018-2141-z

[54]. Brkić D. Solution of the implicit Colebrook equation for flow friction using Excel. Spreadsheets Educ. (EJSiE), 10:4663, 2017. https://sie.scholasticahq.com/article/4663

[55]. Botchkarev A. Assessing Excel VBA suitability for Monte Carlo simulation. Spreadsheets Educ. (EJSiE), 8:4629, 2015. https://sie.scholasticahq.com/article/4629

[56]. Christie D.S. Build your own Monte Carlo spreadsheet. Spreadsheets Educ. (EJSiE), 11:4672, 2018. https://sie.scholasticahq.com/article/4672

[57]. Cordero A., Martí P., Victoria M. Optimización de topología robusta de estructuras continuas usando el método de Monte Carlo y modelos Kriging (Robust topology optimization of continuum structures using Monte Carlo method and Kriging models) /in Spanish/. Rev. Int. Métodos Numér. Cálc. Diseño Ing., 34:5, 2018. doi:10.23967/j.rimni.2017.5.005

[58]. Sobol I.M., Turchaninov V.I., Levitan Y.L., Shukhman B.V. Quasi-random sequence generators; Distributed by OECD/NEA data bank; Keldysh Institute of Applied Mathematics; Russian Academy of Sciences: Moscow, Russia. 1992. https://ec.europa.eu/jrc/sites/jrcsh/files/LPTAU51.rar

[59]. Praks P., Brkić D. Choosing the optimal multi-point iterative method for the Colebrook flow friction equation. Processes, 6:130, 2018. doi:10.3390/pr6080130

[60]. Praks P., Brkić D. Advanced iterative procedures for solving the implicit Colebrook equation for fluid flow friction. Adv. Civ. Eng., 2018:5451034, 2018. doi:10.1155/2018/5451034

[61]. Sun S., Ouyang R., Zhang B., Zhang T.Y. Data-driven discovery of formulas by symbolic regression. MRS Bull., 44:559-564, 2019. doi:10.1557/mrs.2019.156

[62]. Udrescu S.-M., Tegmark M. AI Feynman: A physics-inspired method for symbolic regression. Sci. Adv. 6:eaay2631, 2020. doi:10.1126/sciadv.aay2631

[63]. Gholamy A., Kreinovich V. How to use absolute-error-minimizing software to minimize relative error: Practitioner's guide. Int. Math. Forum, 12:763-770, 2017. doi:10.12988/imf.2017.7761

[64]. Vatankhah A.R. Approximate analytical solutions for the Colebrook equation. J. Hydraul. Eng., 144:06018007, 2018. doi:10.1061/(ASCE)HY.1943-7900.0001454

[65]. Brkić D., Praks P. Discussion of "Approximate analytical solutions for the Colebrook equation.", by Ali R. Vatankhah. J. Hydraul. Eng., 146:07019011, 2020. doi:10.1061/(ASCE)HY.1943-7900.0001667

[66]. Lamri A.A. Discussion of "Approximate analytical solutions for the Colebrook equation.", by Ali R. Vatankhah. J. Hydraul. Eng., 146:07019012, 2020. doi:10.1061/(ASCE)HY.1943-7900.0001668

[67]. Clamond D. Efficient resolution of the Colebrook equation. Ind. Eng. Chem. Res., 48:3665-3671, 2009. doi:10.1021/ie801626g

[68]. Niazkar M. Discussion of "Accurate and efficient explicit approximations of the Colebrook flow friction equation based on the Wright ω-function" by Dejan Brkić and Pavel Praks. Mathematics, 8, 2020.

[69]. Praks P., Brkić D. Accurate and efficient explicit approximations of the Colebrook flow friction equation based on the Wright ω-function: Reply to the Discussion by Majid Niazkar. Mathematics, 8, 2020.

[70]. Brkić D., Praks P. Colebrook's flow friction explicit approximations based on fixed-point iterative cycles and symbolic regression. Computation, 7:48, 2019. doi:10.3390/computation7030048

[71]. Praks P., Brkić D. Rational Approximation for solving an implicitly given Colebrook flow friction equation. Mathematics, 8:26, 2020. doi:10.3390/math8010026

[72]. Praks P., Brkić D. One-log call iterative solution of the Colebrook equation for flow friction based on Padé polynomials. Energies, 11:1825, 2018. doi:10.3390/en11071825